\documentclass[12pt,reqno]{amsart}
\usepackage{amsmath,amssymb,amsfonts}
\usepackage{eucal,enumerate,graphicx,psfrag}

\numberwithin{equation}{section}

\newtheorem{lem}{Lemma}[section]
\newtheorem{cor}[lem]{Corollary}

\newtheorem{thm}[lem]{Theorem}

\theoremstyle{remark}
\newtheorem{exam}[lem]{Example}
\newtheorem{prob}[lem]{Problem}

\newcommand\ts{\textstyle}
\renewcommand{\phi}{\varphi}
\renewcommand{\epsilon}{\varepsilon}
\newcommand\eset{\varnothing}
\newcommand\setm{\setminus}

\newcommand\inv{^{-1}}
\newcommand\textb{\text{\rm b}}
\newcommand\Lat{\operatorname{Lat}}
\newcommand\Latb{\operatorname{\Lat^{\textb}}}
\newcommand\aff{\operatorname{aff}}
\newcommand\conv{\operatorname{conv}}
\newcommand\rk{\operatorname{rk}}
\newcommand\codim{\operatorname{codim}}
\newcommand\Image{\operatorname{Im}}
\newcommand\vol{\operatorname{vol}}
\newcommand\PH{{P,\cH}}	
\newcommand\PoH{{P^\circ,\cH}}	
\newcommand \full{^{{}^{{}_{{}_\bullet}}}}

\newcommand\cA{\mathcal{A}}
\newcommand\cF{\mathcal{F}}
\newcommand\cH{\mathcal{H}}
\newcommand\cL{\mathcal{L}}
\newcommand\cR{\mathcal{R}}
\newcommand\cS{\mathcal{S}}

\newcommand\bbP{\mathbb{P}}
\newcommand\bbR{\mathbb{R}}
\newcommand\bbZ{\mathbb{Z}}

\newcommand\bE{\mathbf{E}}
\newcommand\bj{\mathbf{1}}
\newcommand\0{\hat 0}
\newcommand\1{\hat 1}

\renewcommand\qedsymbol{\ensuremath{\blacksquare}}

\begin{document}

\title{Inside-Out Polytopes}

\author{Matthias Beck}
\author{Thomas Zaslavsky}

\begin{abstract}
 We present a common generalization of counting lattice points in rational 
polytopes and the enumeration of proper graph colorings, nowhere-zero 
flows on graphs, magic squares and graphs, antimagic squares and graphs, 
compositions of an integer whose parts are partially distinct, and 
generalized latin squares.  Our method is to generalize Ehrhart's theory 
of lattice-point counting to a convex polytope dissected by a hyperplane 
arrangement.  We particularly develop the applications to graph and 
signed-graph coloring, compositions of an integer, and antimagic 
labellings.
 \end{abstract}



\thanks{Part of the work of the first author was done while he was a
Robert Riley Assistant Professor at Binghamton University, SUNY.  He thanks the Department of Mathematical Sciences at Binghamton for its hospitality.}
\thanks{The research of the second author was partially supported by National Science Foundation grant DMS-0070729.}

\date{Version of 13 August 2005.}

\maketitle

\noindent Matthias Beck,
\newline\noindent Department of Mathematics,
\newline\noindent San Francisco State University,
\newline\noindent 1600 Holloway Avenue,
\newline\noindent San Francisco, CA 94132, U.S.A.
\newline\noindent Email: {\tt beck@math.sfsu.edu}
\bigskip

\noindent Thomas Zaslavsky,
\newline\noindent Department of Mathematical Sciences,
\newline\noindent Binghamton University (SUNY),
\newline\noindent Binghamton, NY 13902-6000, U.S.A. 
\newline\noindent Email: {\tt zaslav@math.binghamton.edu}


\tableofcontents

\section{In which we introduce polytopes, hyperplanes, and lattice points}  \label{intro}

We study lattice-point counting in polytopes with boundary on the inside.  
To say this in a less mysterious way:  we consider a convex polytope, $P$, 
together with an arrangement of hyperplanes, $\cH$, that dissects the 
polytope, and we count points of a discrete lattice, such as the integer 
lattice $\bbZ^d$, that lie interior to $P$ but not in any of the 
hyperplanes.  We refer to the pair $(P,\cH)$ as an \emph{inside-out 
polytope} because the hyperplanes behave like additional boundary inside 
$P$.

\begin{figure}[htbp]
\begin{center}
\includegraphics[totalheight=3in]{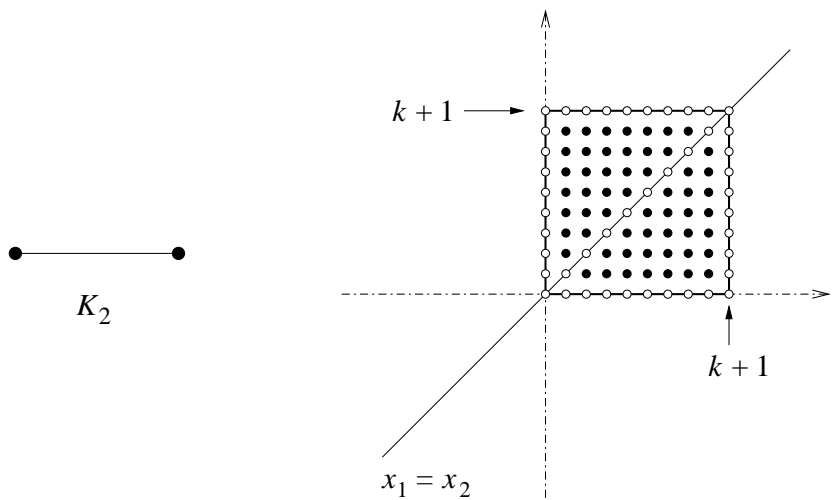} 
\end{center}
\caption{The lattice points in $(k+1) [0,1]^2$ that $k$-color the graph $K_2$, with $k=8$.} \label{F:gcol2}
\end{figure}

We became interested in inside-out lattice-point counting because of a geometrical interpretation of coloring of graphs and signed graphs.
A \emph{coloring} in $c$ colors of a graph $\Gamma$, with node set $V=[n] := \{1,2,\ldots,n\}$, is a function $x: V \to [c]$.  (By $[k]$ we mean the set $\{1,2,\ldots,k\}$, the empty set if $k=0$.)  The coloring $x$ is \emph{proper} if, whenever there is an edge $ij$, $x_i \neq x_j$.  It is a short step to regard $x$ as a point in the real affine space $\bbR^n$ and call it proper if it lies in none of the hyperplanes $h_{ij}: x_i=x_j$ for $ij \in E$, the edge set of $\Gamma$.  That is, if we write 
$$
\cH[\Gamma] := \{ h_{ij}: ij \in E \},
$$ 
 which is the hyperplane arrangement of the graph $\Gamma$, then counting 
proper colorings of $\Gamma$ means counting integral points in $[c]^n 
\setm \bigcup\cH[\Gamma]$, the first instance of an inside-out polytope 
(see Figure \ref{F:gcol2}).  It is well known that the number of proper 
colorings is a polynomial function of $c$, 
$\chi_\Gamma(c)$, called the \emph{chromatic polynomial} of $\Gamma$.
 A famous theorem of Stanley's \cite{AOG} states that when one evaluates 
the chromatic polynomial at negative integers, one obtains the function, \emph{a priori} unrelated to proper graph coloring, that counts pairs consisting of $k$-colorings and compatible acyclic orientations of the graph; in particular, the evaluation at $-1$ gives the number of acyclic orientations.  We will see in Section \ref{graphs} that this fact is a particular case of the general geometrical phenomenon of \emph{Ehrhart reciprocity}, a fundamental theorem in classical lattice-point enumeration in polytopes.

Our purpose in this paper and its sequels \cite{NNZ,MML,SLS} is to apply 
the framework of inside-out polytopes to a multitude of counting problems 
in which there are forbidden values or relationships amongst the values of 
an integral linear function on a finite set which might, for instance, be 
the edge set of a graph or the set of cells of an $n\times n$ square.  
Main examples, aside from graph coloring, are nowhere-zero integral flows, 
magic, antimagic, and latin squares, magic and antimagic graphs, 
compositions (ordered partitions) of an integer into parts with arbitrary 
pairs of parts required to be distinct, and generalizations involving 
rational linear forms.  Our results are of three kinds: 
(quasi)polynomiality of counting functions, M\"obius inversion formulas, 
and the appearance of quantities that generalize the number of acyclic 
orientations of a graph but whose combinatorial interpretation is, in some 
examples at any rate, an open problem.  Among our applications, two stand 
out.  We show how to count antimagic labellings in a systematic way 
(Section \ref{antimagic}), and we explain why a signed graph has not 
one, as with ordinary graphs, but two different chromatic polynomials 
(Theorem \ref{T:sg}).

We have two techniques for attacking the problem of inside-out polytope 
counting.  In the first we dissect the polytope into its intersections 
with the regions of the hyperplane arrangement.  The intersections are 
rational polytopes whose Ehrhart (quasi)polynomials sum to that of the 
inside-out polytope.  Thus we deduce (quasi)polynomiality and reciprocity 
together with interpretations of the leading coefficient and constant 
term.  The second technique is M\"obius inversion over the lattice of 
flats of the arrangement, i.e., sophisticated inclusion-exclusion.  The 
strongest results come in applications where the two methods meld, as most 
neatly in graph coloring.  (Curiously, both techniques were anticipated to 
an extent by Stanley, as we recently learned \cite{StanleyPC}.  Stanley 
used a method equivalent to dissection to give a second proof of his 
combinatorial interpretation of the chromatic polynomial at negative 
arguments; the proof is that via the order polynomial in \cite{AOG}.  Much 
later, Kochol applied a dissection argument to nowhere-zero flows 
\cite{Kochol}.  Then, in his textbook \cite[Exercise 4.10]{EC1} Stanley 
suggests M\"obius inversion over the Boolean algebra or the partition 
lattice to find the number of nonnegative integral solutions $x$, with all 
coordinates distinct, of a rational linear system $Ax=0$---such as the 
equations of a magic square.)

We conclude this paper with two short sections on supplemental topics: subspace arrangements and general valuations.  These are intended to clarify the phenomena by indicating the essential requirements for a theory of our type.  A lattice-point count is one kind of valuation; another example is the \emph{combinatorial Euler characteristic}, which is the alternating sum $a_0 - a_1 + \cdots$ of the number $a_i$ of open cells of each dimension $i$ into which a geometrical object can be decomposed.  In Section \ref{val} we show that the M\"obius inversion formulas (as will be no surprise) apply to any valuation.

\section{In which more characters take the stage}  \label{defs}

We first expand on the geometry of real hyperplane arrangements.  A 
\emph{hyperplane arrangement} $\cH$ is a set of finitely many linear or 
affine hyperplanes in $\bbR^d$.  It divides up the space into regions: an 
\emph{open region} is a connected component of $\bbR^d \setm \bigcup\cH$ 
and a \emph{closed region} is the topological closure of an open region.  
The number of regions into which a hyperplane arrangement $\cH$ divides $\bbR^d$ is $(-1)^d p_\cH(-1)$ \cite{FUTA}, where $p_\cH$ is the characteristic polynomial of $\cH$, defined below.

The \emph{M\"obius function} of a finite partially ordered set (a \emph{poset}) $S$ is the function $\mu: S \times S \to \bbZ$ defined recursively by 
\begin{equation*}
\mu(r,s) := \begin{cases}
	0			&\text{if } r \not\leq s, \\
	1			&\text{if } r = s, \\
	-\sum_{r \leq u < s} \mu(r,u)	&\text{if } r < s.
	\end{cases}
\end{equation*}
Sources are, \emph{inter alia}, \cite{FCT} and \cite[Section 3.7]{EC1}.  $S$ may be the class of closed sets of a closure operator; in that case if $\eset$ is not closed we define $\mu(\eset, s) := 0$ for $s \in S$.  

In a poset $P$, $\0$ denotes the minimum element and $\1$ the maximum 
element, if they exist.  A \emph{lattice poset} (commonly called simply a 
``lattice'' but we must differentiate it from a discrete lattice) is a 
poset in which any two elements have a greatest lower bound (their meet) 
and a least upper bound (or join).  A \emph{meet semilattice} is a poset 
in which meets exist but not necessarily joins.

One kind of poset is the \emph{intersection semilattice} of an affine 
arrangement of hyperplanes, namely,
 $$
\cL(\cH) := \big\{ {\ts\bigcap} \cS : \cS \subseteq \cH \text{ and } {\ts\bigcap} \cS \neq \eset \big\},
$$ 
ordered by reverse inclusion \cite{FUTA}.
The elements of $\cL(\cH)$ are called the \emph{flats} of $\cH$.  $\cL$ is a \emph{geometric semilattice} (of which the theory is developed in \cite{McNaff,WW}) with $\0 = \bigcap \eset = \bbR^d$; it is a \emph{geometric lattice} (for which see \cite[p.\ 357]{FCT} or \cite{EC1}, etc.) if $\cH$ has nonempty intersection, as when all the hyperplanes are homogeneous.  
A hyperplane arrangement decomposes the ambient space into relatively open cells called \emph{open faces} of $\cH$, whose topological closures are the \emph{closed faces}.  For a more precise definition we need the arrangement \emph{induced} by $\cH$ on a flat $s$; this is
$$
\cH^s := \{ h\cap s : h\in\cH, h\not\supseteq s \}.
$$
A face of $\cH$ is then a region of any $\cH^s$ for $s \in \cL(\cH)$.  One face is $\1 = \bigcap \cH$ itself, if nonempty.  
An oddity about hyperplane arrangements is that, for technical reasons, one wants to treat the whole space as a hyperplane (called the \emph{degenerate hyperplane}) that may or may not belong to $\cH$.  If $\cH$ contains the degenerate hyperplane, it has no regions, because $\bbR^d\setm \bigcup\cH = \eset$.  However, $\cH$ does have faces; e.g., its $d$-dimensional faces are the regions of $\cH^{\0}$, the arrangement induced in $\0 = \bbR^d$.  

The \emph{characteristic polynomial} of $\cH$ is defined in terms of the M\"obius function of $\cL(\cH)$ by
$$
p_\cH(\lambda) := \begin{cases}
0	&\text{if $\cH$ contains the degenerate hyperplane, and} \\
\sum_{s \in \cL(\cH)} \mu(\0,s) \lambda^{\dim s}	&\text{otherwise}.
\end{cases}
$$
And three more definitions:  for a set or point $X$ in $\bbR^d$, 
\begin{equation*}
\cH(X) := \{h \in \cH : X\subseteq h\},
\end{equation*}
\begin{equation*}
s(X) := \bigcap \cH(X) = \bigcap \{h \in {\cH} : X\subseteq h\},
\end{equation*}
the smallest flat of $\cH$ that contains $X$, and 
$$
F(X) := \text{the unique open face of $\cH$ that contains $X$}
$$
provided that $X$ is contained in an open face, as for instance when it is a point.

A \emph{convex polytope} $P$ is a bounded, nonempty set that is the intersection of a finite number of open and closed half spaces in $\bbR^d$; $P$ may be closed, relatively open, or neither.  
A closed convex polytope is also the convex hull of a finite set of points in $\bbR^d$.  
Another kind of poset is the face lattice of $P$.  
A \emph{closed face} of $P$ is either $\bar P$, the topological closure of $P$, or the intersection with $\bar P$ of any hyperplane $h$ such that $P \setm h$ is connected.  (Then $h$ is a \emph{supporting hyperplane} of $P$; this includes hyperplanes that do not intersect $P$ at all.)  
An \emph{open face} is the relative interior $F^\circ$ of a closed face $F$.  
(The relative interior of a point is the point.)  
The null set is a face; it and $\bar P$ (or $P^\circ$) are the \emph{improper} faces.  
The \emph{face lattice} $\cF(P)$ is the set of open faces, partially ordered by inclusion of the closures.  
A \emph{vertex} is a 0-dimensional face.  
A \emph{facet} is a face whose dimension is $\dim P-1$; a \emph{facet hyperplane} is the affine span of a facet.  If $P\subseteq \bbR^d$ is not full-dimensional, then a facet hyperplane is a relative hyperplane of the affine flat spanned by $P$.

A \emph{dilation} of a set $X\subseteq \bbR^d$ is any set $tX = \{tx : x\in X\}$ for a real number $t>0$.

\section{In which we encounter facially weighted enumerations}  \label{conv}

Our first main result expresses the Ehrhart quasipolynomials of an inside-out polytope $(P,\cH)$ in terms of the combinatorics of $\cH$; but its natural domain is far more general.  We may take any discrete set $D$ in $\bbR^d$, any bounded convex set $C$, and any hyperplane arrangement $\cH$ that is \emph{transverse} to $C$: every flat $u\in \cL(\cH)$ that intersects the topological closure $\bar C$ also intersects $C^\circ$, the relative interior of $C$, and $C$ does not lie in any of the hyperplanes of $\cH$.  
A convenient sufficient condition for transversality is that $C^\circ \cap \bigcap \cH \neq \eset$ and $C \not\subseteq \bigcup\cH$.

A \emph{region of $(C,\cH)$}, or \emph{of $\cH$ in $C$}, is one of the components of $C \setm \bigcup\cH$, or the closure of such a component.  
A \emph{vertex} of $(C,\cH)$ is a vertex of any of its regions.  

The \emph{multiplicity of $x\in \bbR^d$ with respect to $\cH$} is 
\begin{align*}
m_\cH(x) &:= \text{the number of closed regions of $\cH$ that contain $x$}.
\intertext{The \emph{multiplicity with respect to $(C,\cH)$} is }
m_{C,\cH}(x) &:= \begin{cases}
\text{the number of closed regions of $(C,\cH)$ that contain $x$}, 
  & \text{if } x \in C, \\
0, & \text{if } x \notin C.
\end{cases}
\end{align*}
This may not equal $m_\cH(x)$ for $x \in C$, unless one assumes transversality.  
The \emph{closed} and \emph{open $D$-enumerators} of $(C,\cH)$ are 
\begin{align*}
E_{C,\cH}(D) &:= \sum_{x\in D} m_{C,\cH}(x)
\intertext{and}
E^\circ_{C,\cH}(D) &:= \# \big( D\cap C \setm \ts\bigcup\cH \big).
\end{align*}

\begin{thm}  \label{T:mobius}  
Let $C$ be a full-dimensional, bounded, convex subset of $\bbR^d$, $\cH$ a hyperplane arrangement not containing the degenerate hyperplane, and $D$ a discrete set in $\bbR^d$.  Then
\begin{equation}\label{E:omobius}
E^\circ_{C,\cH}(D) = \sum_{u\in \cL(\cH)}\mu(\0,u)\, \#(D\cap C\cap u)
\end{equation}
and if $\cH$ is transverse to $C$, 
\begin{equation}\label{E:mobius}
E_{C,\cH}(D) = \sum_{u\in \cL(\cH)} |\mu(\0,u)|\, \#(D\cap C\cap u),
\end{equation}
where $\mu$ is the M\"obius function of $\cL(\cH)$.  
\end{thm}

Note that $u$ can be omitted from the sum if $u\cap C=\eset$.  We formalize this by defining 
$$
\cL(C,\cH) := \{ u \in \cL(\cH) : u\cap C \neq \eset \} ,
$$
 the \emph{intersection poset} of $(C,\cH)$, and observing that 
$\cL(C,\cH)$ can replace $\cL(\cH)$ in the range of summation of either 
equation and $\cL(C^\circ,\cH)$ can replace $\cL(\cH)$ in Equation 
\eqref{E:mobius}.  This is helpful in solving examples.

\begin{proof}[Proof of Equation \eqref{E:omobius}]
We begin with the observation that, for any flat $r$ of $\cH$, 
\begin{equation*}
\#(D\cap C\cap r) = \sum_{u \in \cL: u\geq r} E^\circ_{C\cap u,\cH^u}(D).
\end{equation*}
The reason for this is that $C\cap r$ is the disjoint union of all the open faces of $\cH$ in $C\cap r$, $C\cap u \setm \bigcup\cH^u$ is the disjoint union of the open faces of $\cH$ that span $u$, and counting points of $D$ is an additive function on open faces (a \emph{valuation} in technical language).  By M{\"o}bius inversion,
\begin{equation*}
E^\circ_{C\cap r,\cH^r}(D)  = \sum_{u\in \cL: u\geq r} \mu(r,u) \#(D\cap C\cap S).
\end{equation*}
Setting $r = \0$ gives the desired formula unless $\cH$ contains the degenerate hyperplane.  In that case, however, both sides of \eqref{E:omobius} equal zero.
\end{proof}

The proof of the second equation depends on two lemmas about transversality and an algebraic expression for the multiplicity of a point.  

\begin{lem}  \label{L:flats} 
Let $C$ be a convex set and $\cH$ a transverse hyperplane arrangement.  If $u$ is a flat of $\cH$, then $\overline{u\cap C^\circ} = u\cap \bar C$.  
\end{lem}

\begin{proof}  
We need to prove that every neighborhood of a point $x\in u\cap \partial C$ intersects $u\cap C^\circ$.  By transversality $u$ intersects $C^\circ$, say in a point $y$.  Then the segment $\conv(x,y)$ lies in $u\cap C^\circ$ except possibly for $x$.  This implies our desideratum. 
\end{proof}

\begin{lem}  \label{L:faces}  
Let $C$ be a convex set and $\cH$ a transverse hyperplane arrangement.  If $F$ is a face of $\cH$, then $\overline{F\cap C^\circ} = \bar F\cap \bar C$.  
\end{lem}  

\begin{proof}  
The question reduces to proving that $\overline{F\cap C^\circ} \supseteq \bar F\cap \bar C$ when $F$ is an open face of $\cH$.  Take $x\in (\bar F\cap \bar C)\setm (F\cap C^\circ)$.  Then $x\in F(x)\cap \partial C$.  
By Lemma \ref{L:flats} with $u = \aff F(x)$, every neighborhood of $x$ intersects $u\cap C^\circ$.  Because $F(x)$ is open in $u$, every small neighborhood of $x$ in $u$ is contained in $F(x)$.  Since $F(x) \subseteq \bar F$, every small neighborhood of $x$ in $\bbR^d$ meets $F\cap C^\circ$ and therefore $x\in \overline{F\cap C^\circ}$.  
\end{proof}  

\begin{lem} \label{L:hmult}
The multiplicity of $x\in \bbR^d$ with respect to a hyperplane arrangement $\cH$ is given by
\begin{equation*}
m_{\cH}(x) = (-1)^{\codim s(x)}p_{\cH(x)}(-1).
\end{equation*}
\end{lem}

\begin{proof}  $x$ belongs to the unique open face $F(x)$, which is an open region of $\cH^{s(x)}$. There is an obvious bijection between (closed) regions $R$ of $\cH$ that contain $F(x)$ and regions $R'$ of $\cH(x)$:  $R' \leftrightarrow R$ if $R' \subseteq R$. In fact, each region $R'$ contains $s(x)$ and is dissected by $\cH\setm \cH(x)$ into regions of $\cH$, of which one and only one contains $x$.  Therefore, the number of regions of $\cH$ that contain $x$ equals the number of regions of $\cH(x)$, which is $(-1)^{\codim s(x)} p_{\cH(x)}(-1)$.  
\end{proof}

\begin{lem} \label{L:mult}
Let $C$ be a full-dimensional, bounded, convex subset of $\bbR^d$ and $\cH$ a transverse hyperplane arrangement. The multiplicity of $x\in \bbR^d$ with respect to $C$ and $\cH$ is given by
\begin{equation*}
m_{C,\cH}(x) = \begin{cases} 
(-1)^{\codim s(x)}p_{\cH(x)}(-1) &\text{if $x\in C$,} \\
0 &\text{if $x\not\in C$.}
\end{cases}
\end{equation*}
\end{lem}

\begin{proof}  
We may assume $x\in C$.  In relation to $C$, $x$ lies in the open face $F(x)\cap C$ of $\cH$ in $C$.  
We must prove that every region $R$ of $\cH$ that contains $x$ corresponds to a region of $\cH$ in $C$, so that $m_{C,\cH}(x)=m_\cH(x)$.  
This follows from Lemma \ref{L:faces} with $F = R^\circ$:  since $x \in R \cap \bar C$, $R^\circ \cap C^\circ$ is nonempty, so $R\cap C$ is full-dimensional.  
\end{proof}

\begin{proof}[Proof of Equation \eqref{E:mobius}] 
We apply Lemma \ref{L:mult}, the definition of the characteristic polynomial, and Rota's sign theorem \cite[Section 7, Theorem 4]{FCT}.
\begin{equation*}
\begin{aligned}
E_{C,\cH}(D) &= \sum_{x\in C\cap D} |p_{\cH(x)}(-1)|\\
&= \sum_{x\in C\cap D}\ \sum_{u\leq s(x)} |\mu(\0,u)|\\
&= \sum_{u\in \cL} |\mu(\0,u)| \#\{x\in C\cap D: s(x) \subseteq u\}\\
&= \sum_{u\in \cL}|\mu(\0,u)| \# (D\cap C\cap u).
\end{aligned}
\end{equation*}
\end{proof}

There is also a proof of \eqref{E:mobius} by inversion.  See Section \ref{val}.

\section{In which we arrange Ehrhart theory with hyperplanes}  \label{ehr}

A \emph{discrete lattice} is a set of points in an affine space (such as $\bbZ^d$ in $\bbR^d$) that is locally finite and is invariant under any translation that carries some lattice point to another lattice point.  
We call a polytope \emph{$D$-integral} if its vertices all lie in $D$ and \emph{$D$-fractional} if the vertices lie in the contracted lattice $t\inv D$ for some positive integer $t$.  
The \emph{denominator} of $P$ \cite{Ebook} is the smallest such $t$.  
(To define $t\inv D$ here we assume coordinates chosen so that $0 \in D$. Later, at Corollary \ref{C:affine}, we deal with a more general situation.)  
A \emph{quasipolynomial} is a function $q(t) = \sum_0^d c_it^i$ with coefficients $c_i$ that, though not necessarily constant, at any rate are periodic functions of $t$ (so that $q(t)$ is a polynomial on each residue class modulo some integer, called the \emph{period}; these polynomials are the \emph{constituents} of $q$). 
According to Ehrhart \cite{EhrCR,Ehr}, if $P$ is a closed, $d$-dimensional, $D$-fractional convex polytope and 
\begin{equation}\label{E:ehrdef}
E_P(t) := \#(D \cap tP) = \#(P \cap t\inv D),
\end{equation}
 then $E_P$ is a quasipolynomial whose degree is $d$, whose period divides 
the denominator of $P$, and whose leading coefficient equals $\vol_D {P}$, 
the volume of $P$ normalized with respect to $D$ (that is, we take the 
volume of a fundamental domain of $D$ to be 1; in the case of the integer 
lattice $\bbZ^d$ this is the ordinary volume), and whose \emph{constant 
term} $E_P(0)$ equals $1$ \cite{Ehr,MacdV}.  Defining the Ehrhart 
quasipolynomial of any rational polytope $P$, not necessarily closed, by 
Equation \eqref{E:ehrdef}, it applies to relatively open as well as closed 
polytopes, except that then the constant term is the combinatorial Euler 
characteristic of $P$.  Ehrhart \cite{Ehr} then conjectured and he, 
Macdonald \cite{MacdP}, and McMullen \cite{McMul} proved the 
\emph{reciprocity law}
 $$
 E_{P^\circ}(t) = (-1)^{\dim P}E_{P}(-t). 	\label{E:ehrrecip}
$$

Our theory begins with a rational, closed convex polytope $P$ and an arrangement $\cH$ of rational hyperplanes that is transverse to $P$.  Rationality means that the vertices of $P$ are rational points and the hyperplanes in $\cH$ are specified by equations with rational coefficients.  We call $(P,\cH)$ a \emph{rational inside-out polytope} of \emph{dimension} $\dim P$.  
More generally we have any discrete lattice $D$, a closed $D$-fractional convex polytope $P$, and a $D$-fractional hyperplane arrangement $\cH$ (transverse to $P$): that is, each hyperplane in $\cH$ is spanned by the $D$-fractional points it contains.  Then $(P,\cH)$ is a \emph{$D$-fractional inside-out polytope}.  
The \emph{vertices of $(P,\cH)$} are all the intersection points in $P$ formed by the hyperplanes of $\cH$ and the facets of $P$, including the vertices of $P$. 
The \emph{denominator} of $(P,\cH)$ (with respect to the discrete lattice $D$) is the smallest positive integer $t$ for which $t\inv D$ contains every vertex of $(P,\cH)$.
We call $(P,\cH)$ \emph{$D$-integral} if all its vertices lie in $D$.  
We always assume that $P$ is closed.

The Ehrhart quasipolynomials of $(P,\cH)$ are the (\emph{closed}) \emph{Ehrhart quasipolynomial}, 
\begin{align*}
E_{P,\cH}(t) &:= 
\sum_{x\in t\inv D}m_{P,\cH}(x),
\intertext{and the \emph{open Ehrhart quasipolynomial}, }
E^\circ_{P,\cH}(t) &:= 
\# \big( t\inv D \cap \big[ P \setm \ts\bigcup\cH \big] \big),
\end{align*}
both defined for positive integers $t$ in terms of the $D$-enumerators of Section \ref{conv}, with $D$ replaced by $t\inv D$.  Thus if $P$ is full-dimensional and $R_1,\hdots,R_k$ are the closed regions of $(P,\cH)$,
\begin{equation} \label{E:sum}
E_{P,\cH}(t) = \sum^k_{i=1} E_{R_i}(t)\quad \text{ and }\quad E^\circ_{P^\circ,\cH}(t) = \sum^k_{i=1} E_{R^\circ_i}(t) \ .
\end{equation}

\begin{thm}\label{T:quasi}
If $D$ is a full-dimensional discrete lattice and $(P,\cH)$ is a closed, full-dimensional, $D$-fractional inside-out polytope in $\bbR^d$ such that $\cH$ does not contain the degenerate hyperplane, then $E_{P,\cH}(t)$ and $E^\circ_{P^\circ,\cH}(t)$ are quasipolynomials in $t$, with period equal to a divisor of the denominator of $(P,\cH)$, with leading term $c_d t^d$ where $c_d = \vol_D {P}$, and with the constant term $E_{P,\cH}(0)$ equal to the number of regions of $(P,\cH)$.  Furthermore,
\begin{equation}\label{E:reciprocity}
E^\circ_{P^\circ,\cH}(t) = (-1)^d E_{P,\cH}(-t).
\end{equation}
\end{thm}

\begin{proof} 
By \eqref{E:sum}, standard Ehrhart theory, and the fact that a closed region has Euler characteristic 1.
\end{proof}

The periodically varying \emph{quasiconstant term} $c_0(t)$ has no presently known interpretation, save at $t\equiv0$.

One might use Theorem \ref{T:quasi} to compute the number of regions of a hyperplane arrangement.  Suppose, for instance, that $\cH$ has nonempty intersection and this intersection meets the interior of $P$.  Then the constant term $E_{P,\cH}(0)$ equals the number of regions of $\cH$.  If one can evaluate $E_{P,\cH}(kp)$ (where $p$ is the period) for enough values of $k$, one can deduce the constant term, thus the number of regions, by polynomial interpolation.  Sometimes this is feasible, as with the simpler examples in \cite{SLS}.

It is easy to prove as well that $E^\circ_{P,\cH}(t)$ and $E_{P^\circ,\cH}(t)$ are quasipolynomials in $t$ with some of the same properties, e.g., the leading term, although we do not know they have the same period as each other or as $E_{P,\cH}(t)$.

The first theorem does not require transversality, but the next one does, in part.

\begin{thm}\label{T:ehrhyp} If $D$, $P$, and $\cH$ are as in Theorem \ref{T:quasi}, except that $P$ need not be closed, then
\begin{equation}\label{E:oehrhyp}
E^\circ_{\PoH}(t) = \sum_{u\in \cL(\cH)} \mu(\0,u) E_{P^\circ\cap u}(t) ,
\end{equation}
and if $\cH$ is transverse to $P$,
\begin{equation} \label{E:ehrhyp}
E_{\PH}(t) = \sum_{u\in \cL(\cH)}|\mu(\0,u)|E_{P\cap u}(t) .
\end{equation}
\end{thm}

\begin{proof} 
A special case of Theorem \ref{T:mobius}.
\end{proof}

The range of summation may be taken to be the intersection poset $\cL(\PoH)$ if one prefers a smaller sum.

If it so happens that, as in the graph coloring examples, $E_{P\cap u}(t)
= f(t)^{\dim u}$, then the right side of \eqref{E:ehrhyp} becomes $(-1)^d       
p_{\cH}(-f(t))$ and that of \eqref{E:oehrhyp} becomes $p_{\cH}(f(t))$.

 When computing specific examples (as in \cite{SLS}) we find it most 
convenient to work with generating functions; thus we need the generating 
function version of Theorem \ref{T:mobius}.  Define
\begin{align*}
\bE^\circ_{C,\cH}(x) = \sum_{t=1}^\infty E^\circ_{C,\cH}(t)\, x^t , \qquad
\bE_{C,\cH}(x) = \sum_{t=0}^\infty E_{C,\cH}(t)\, x^t ,
\end{align*}
 where $C$ is a closed or relatively open convex polytope.  Ehrhart 
reciprocity (Theorem \ref{T:quasi}) is expressed as
 \begin{equation}\label{E:recipgf}
 \bE^\circ_{\PoH}(x) = (-1)^{1+\dim P}\, \bE_{\PH}(x\inv) ,
 \end{equation}
 proved by summing over all regions of $(\PH)$ the ordinary 
generating-function reciprocity formula
 $$ 
 \bE^\circ_{P^\circ}(x) = (-1)^{1+\dim P}\, \bE_{P}(x\inv)
 $$
 (see \cite[Theorem 4.6.14]{EC1}).  M\"obius summation (Theorem 
\ref{T:ehrhyp}) becomes
\begin{equation}\label{E:oehrhypgf}
\bE^\circ_{\PoH}(t) = 
\sum_{u\in \cL(\cH)} \mu(\0,u)\, \bE^\circ_{P^\circ\cap u}(t) ,
\end{equation}
and, if $\cH$ is transverse to $P$,
\begin{equation} \label{E:ehrhypgf}
\bE_{\PH}(t) = 
\sum_{u\in \cL(\cH)} |\mu(\0,u)|\, \bE_{P\cap u}(t) .
\end{equation}
 (As in Theorem \ref{T:ehrhyp}, the range of summation may be taken to 
be $\cL(\PoH)$ if preferred.) 
 These two equations are immediate from Theorem \ref{T:ehrhyp} except 
for the constant term of \eqref{E:ehrhypgf}.  By transversality we may 
replace $\cL(\cH)$ by $\cL(\PoH)$.  Quasipolynomiality implies that 
Equation \eqref{E:ehrhyp} holds for all integers $t$, in particular $t=0$.

One can potentially evaluate the inside-out Ehrhart quasipolynomial $E^\circ_\PoH$ in an example by counting lattice points for enough values of $t$ and using polynomial interpolation.  (This requires knowing an upper bound on the period, such as the denominator.)  The number of values necessary is lessened if one knows something of the coefficients in advance.  
For example, the leading coefficient of every constituent of $E^\circ_\PoH$ is $\vol P$; this is also the leading coefficient of the ordinary Ehrhart quasipolynomial $E_{P^\circ}$, which can be interpolated from many fewer evaluations.  
By Theorem \ref{T:ehrhyp} one can simplify the computation further, even without evaluating the M\"obius function, if one first calculates $E_{P^\circ}$. Consider the second leading coefficients $c_{d-1}(t)$ in $E^\circ_\PoH$ and $c_{P,d-1}(t)$ in $E_{P^\circ}$.  
By Theorem \ref{T:ehrhyp},
\begin{equation} \label{E:subleading}
c_{d-1}(t) = c_{P,d-1}(t) - \sum_{\substack{u\in \cL(\PoH)\\ \codim u = 1}} \vol(P\cap u) ,
\end{equation}
because $\mu(\0,u)=-1$.  The sum is a constant, so if it is evaluated for one constituent of $E^\circ_\PoH$ it is known for all; whence one needs fewer evaluations to determine the coefficients of all the constituents of $E^\circ_\PoH$.  
This idea sees a practical application in one of the methods in \cite{SLS}.

Sometimes (as in \cite{MML,SLS}) the polytope is not full-dimensional; its affine span, $\aff P$, 
might not even intersect the discrete lattice.  Suppose, then, that $D$ is 
a discrete lattice in $\bbR^d$ and $s$ is any affine subspace.  The 
\emph{period} $p(s)$ of $s$ with respect to $D$ is the smallest positive 
integer $p$ for which $p\inv D$ meets $s$.  Then Theorem \ref{T:quasi} 
implies the following.

\begin{cor} \label{C:affine}
 Let $D$ be a discrete lattice in $\bbR^d$, $P$ a $D$-fractional convex 
polytope, and $\cH$ a hyperplane arrangement in $s := \aff P$ that does not contain the degenerate hyperplane.  
Then $E_{P,\cH}(t)$ and $E^\circ_{P^\circ,\cH}(t)$ are quasipolynomials in $t$ that satisfy the reciprocity law $E^\circ_{P^\circ,\cH}(t) = (-1)^{\dim s}E_{P,\cH}(-t)$. 
Their period is a multiple of $p(s)$ and a divisor of the denominator of $(P,\cH)$.  
If $t \equiv 0 \mod {p(s)}$, the leading term of $E_{P,\cH}(t)$ is $(\vol_{p(s)\inv D}{P}) t^{\dim s}$ and its constant term is the number of regions of $(P,\cH)$; but if $t \not\equiv 0 \mod {p(s)}$, then $E_{P,\cH}(t) = E^\circ_{P^\circ,\cH}(t) = 0$.
\hfill$\qedsymbol$
 \end{cor}

The period's being greater than one suggests that we should renormalize,  multiplying $s$, $P$, and $\cH$ by $p(s)$.  This divides both the denominator of the inside-out polytope 
and the period of the Ehrhart quasipolynomials by $p(s)$ and eliminates the zero constituents of the quasipolynomials.  The $3\times3$ magic squares are a perfect example \cite{SLS}.

Usually $\cH$ will be induced by an arrangement $\cH_0$ in $\bbR^d$.  It is easy to see that $\cH$ is transverse to $P$ if and only if $\cH_0$ is.

\section{In which we color graphs and signed graphs}  \label{graphs}

\subsection*{Unsigned graphs}
 We begin by deriving from inside-out Ehrhart theory some known results on 
the chromatic polynomial of a graph.  
 An \emph{ordinary graph} is a graph $\Gamma$ whose edges are links (with 
two distinct endpoints) and loops (with two coinciding endpoints); 
multiple edges are permitted.  We treat, always, only finite graphs.  The 
\emph{order} is the number of nodes; we write $n$ for the order of 
$\Gamma$.

\begin{thm}  \label{T:g}  
Let $\Gamma$ be an ordinary graph and let $P = [0,1]^n$.  The closed and open Ehrhart quasipolynomials of $(P,\cH[\Gamma])$ satisfy 
\begin{equation*}
(-1)^n E_{P,\cH[\Gamma]}(-t) = E^\circ_{P^\circ,\cH[\Gamma]}(t) = \chi_\Gamma(t-1).
\end{equation*}
\end{thm}

\begin{proof}  
In $t\inv\bbZ^n$ the points that are counted by $E^\circ_{P^\circ,\cH[\Gamma]}(t)$ are those of $(t\inv\{1,2,\hdots,t-1\})^n$ that do not lie in any forbidden hyperplane.  The number of such points is the number of proper $(t-1)$-colorings of $\Gamma$.  
\end{proof}

\begin{cor}[Birkhoff \cite{Birk} for maps, Whitney \cite{Logical} for graphs] \label{C:gpoly}  
For an ordinary graph $\Gamma$ with no loops, $\chi_\Gamma$ is a monic polynomial of degree $n$.  If $\Gamma$ has a loop, $\chi_\Gamma = 0$.  
\end{cor}

\begin{proof}  
Since $P$ is full-dimensional in $\bbR^n$ and has volume 1, the leading term of $\chi_\Gamma$ is $1x^n$ by Ehrhart theory.  It remains to prove that $(P,\cH[\Gamma])$ has denominator 1, or in other words that $(P,\cH[\Gamma])$ has integer vertices.  This is the next lemma. 
\end{proof}

\begin{lem} \label{L:gdenominator} 
If $\Gamma$ is an ordinary graph, $(P,\cH[\Gamma])$ has integer vertices.
\end{lem}

\begin{proof}  
Because, as is well known, $\cH[\Gamma]$ is a hyperplanar representation of the graphic matroid $G(\Gamma)$, the flats of $\cH[\Gamma]$ correspond to closed subgraphs of $\Gamma$, i.e., to partitions $\pi$ of $V$ into blocks that induce connected subgraphs.  The flat $s(\pi)$ corresponding to $\pi$ is described by $x_i = x_j$ if $i$ and $j$ belong to the same block of $\pi$ (we write $i \underset \pi \sim j$).  A vertex of $(P,\cH[\Gamma])$ is determined by a flat $s(\pi)$ of dimension $k$, say, together with $k$ facet hyperplanes of $P$ that have the form $x_i = a_i \in \{0,1\}$. Obviously, such points are integral.  
\end{proof}

\begin{thm} \label{T:gcharacteristic}
For an ordinary graph $\Gamma$, $\chi_\Gamma(c) = p_{\cH[\Gamma]}(c)$.
\end{thm}

\begin{proof}
We apply Theorem \ref{T:ehrhyp}.  It is easy to see that $E_{P^\circ\cap u}(t) = (t-1)^{\dim u}$.  Therefore, $E^\circ_{P^\circ,\cH[\Gamma]}(t) = p_{\cH[\Gamma]}(t-1)$.  This equals $\chi_\Gamma(t-1)$ by Theorem \ref{T:g}.
\end{proof}

Given an orientation $\alpha$ of $\Gamma$ and a $c$-coloring $x : V\to [c]$, Stanley calls them \emph{compatible} if $x_j\geq x_i$ whenever there is a $\Gamma$-edge oriented from $i$ to $j$, and \emph{proper} if $x_j > x_i$ under the same condition \cite{AOG}.  An orientation is \emph{acyclic} if it has no directed cycles.  From Theorem \ref{T:g} we derive a more unified version of Stanley's second proof of his famous result.

\begin{cor}[Stanley \cite{AOG}] \label{C:gcompatible}  
The number of pairs $(\alpha,x)$ consisting of an acyclic orientation of an ordinary graph $\Gamma$ and a compatible $c$-coloring equals $(-1)^n \chi_\Gamma(-c)$.  
In particular, $(-1)^n \chi_\Gamma(-1) =$ the number of acyclic orientations of $\Gamma$.
\end{cor}

\begin{proof}  
From Theorem \ref{T:g}, 
\begin{equation*}
E_{P,\cH[\Gamma]} (t) = (-1)^n\chi_\Gamma(-(t+1)).
\end{equation*}
What $E_{P,\cH[\Gamma]}(t)$ counts is the number of pairs $(x,R)$ where $x$ is a coloring with color set $\{0,1,\hdots,t\}$, $R$ is a closed region of $\cH[\Gamma]$, and $x\in R$.  Greene observed that regions $R$ correspond with acyclic orientations $\alpha$ in the following way:  $R^\circ$ is determined by converting each equation $x_i = x_j$ corresponding to an edge of $\Gamma$ into an inequality $x_i < x_j$; then in $\alpha$ the edge $ij$ is directed from $i$ to $j$.  (See \cite{AO} or \cite[Section 7]{IWN}.)  The orientation is acyclic because $R^\circ \neq \eset$.  Thus $x$ is compatible with $\alpha$ if and only if $x\in R$.  
The final assertion is an instance of the evaluation $E(0)$ in Theorem \ref{T:quasi}.
\end{proof}

This proof generalizes Greene's geometrical approach to counting acyclic orientations (that is, the case $c=1$); see \cite{AO} or \cite[Section 7]{IWN}.


\subsection*{Signed graphs}
\begin{figure}[h]
\begin{center}
\psfrag{pmk2}[r]{$\pm K_2$}
\includegraphics[totalheight=3.5in]{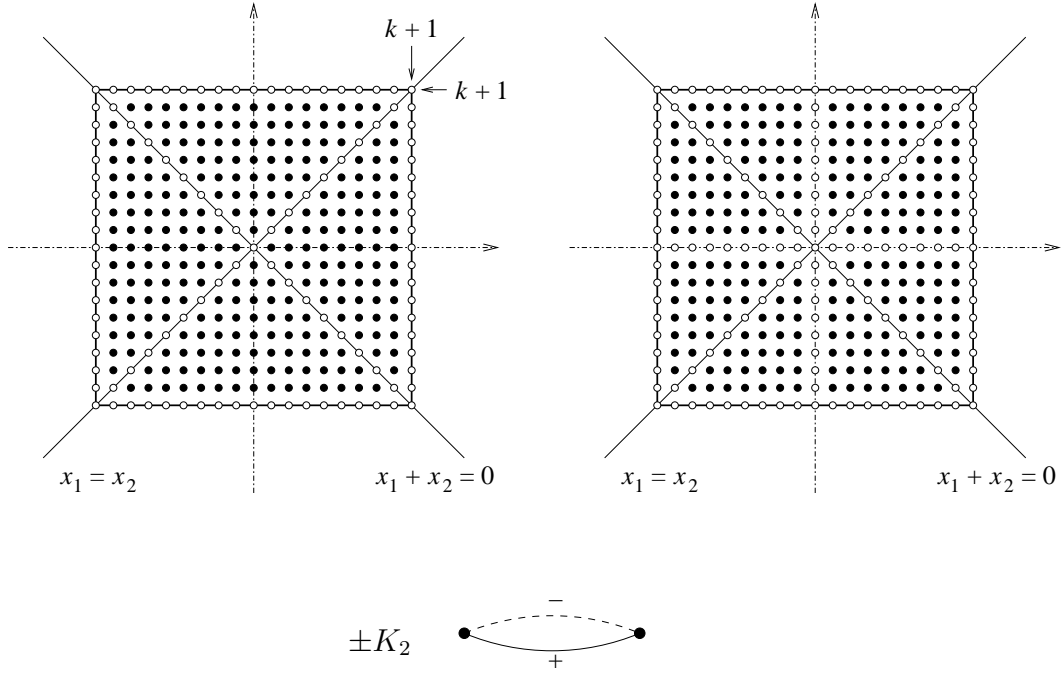}
\end{center}
\caption{Illustrating for $k=8$:  On the left, the lattice points in $(k+1)[-1,1]^2$ that $k$-color the signed graph $\pm K_2$.
On the right, the lattice points that $k$-color it without $0$.} \label{F:sgcol2both}
\end{figure}

A \emph{signed graph} $\Sigma = (\Gamma,\sigma)$ consists of a graph $\Gamma$ (multiple edges allowed) which may have, besides links and loops, also \emph{halfedges} (with only one endpoint) and \emph{loose edges} (no endpoints), and a signature $\sigma$ that labels each link and loop with a sign, $+$ or $-$.  The \emph{order} of $\Sigma$ is the number of nodes, written $n$.
A \emph{$c$-coloring} \cite{SGC} of a signed graph with node set $V=[n]$ is a function 
$$
x: V \to \{ -c,-(c-1),\ldots,0,\ldots,c-1,c \} ;
$$ 
 we say $x$ is \emph{proper} if, whenever there is an edge $ij$ with sign $\epsilon$, then $x_j \neq \epsilon x_i$. 
 Geometrically, $x \in \{ -c,-(c-1),\ldots,c \}^n \setm \bigcup\cH[\Sigma]$ where
 \begin{align*}
\cH[\Sigma] :=\ &\{ h_{ij}^\epsilon : \Sigma \text{ has an edge $ij$ with sign } \epsilon \}\\
 &\cup \{ x_i = 0 : \Sigma \text{ has a halfedge at node } v_i \}\\
 &\cup \{ 0=0 \text{ if $\Sigma$ has a loose edge} \}
 \end{align*}
and $h_{ij}^\epsilon$ is the hyperplane $x_j = \epsilon x_i$.  
(The degenerate hyperplane $0 = 0$ is the set $\bbR^n$, the same as $h_{ii}^+$ belonging to a positive loop.)  
The function 
$$
\chi_\Sigma(2c+1) := \text{the number of proper $c$-colorings of $\Sigma$}
$$ 
is known by \cite{SGC} to be a polynomial, called the \emph{chromatic polynomial} of $\Sigma$; here we prove this from Ehrhart theory.  We see that $\chi_\Sigma(2c+1) = E^\circ_{P^\circ,\cH[\Sigma]}(c+1)$, the number of lattice points that 
lie in $(c+1)P^\circ$ (where now $P=[-1,1]^n$) but not in any of the hyperplanes of $\cH[\Sigma]$.  (See Figure \ref{F:sgcol2both}.)
Furthermore, the regions of $\cH[\Sigma]$ are known to correspond to the acyclic orientations of $\Sigma$ \cite{OSG} and the regions that contain a coloring $x$ correspond to the acyclic orientations that are compatible with $x$ \cite{SGC}.  Thus $E_{P,\cH[\Sigma]}(c)$ is the number of pairs consisting of a coloring and a compatible acyclic orientation, which is known to equal $(-1)^n \chi_\Sigma(-(2c+1))$ \cite{SGC}.

Signed graphs have a second chromatic counting function: the \emph{zero-free chromatic polynomial} 
$$
\chi^*_\Sigma(2c) := \text{the number of proper $c$-colorings } x: V \to \pm[c],
$$ 
that is, it counts colorings not taking the value $0$.  This is obviously also an inside-out Ehrhart polynomial, but it is not obvious that $\chi_\Sigma$ and $\chi^*_\Sigma$ are closely related.  In fact, they are the two constituent polynomials of a single Ehrhart quasipolynomial.

\begin{figure}[h]
\begin{center}
\psfrag{pmk2o}[r]{$\pm K_2^\circ$}
\includegraphics[totalheight=3.5in]{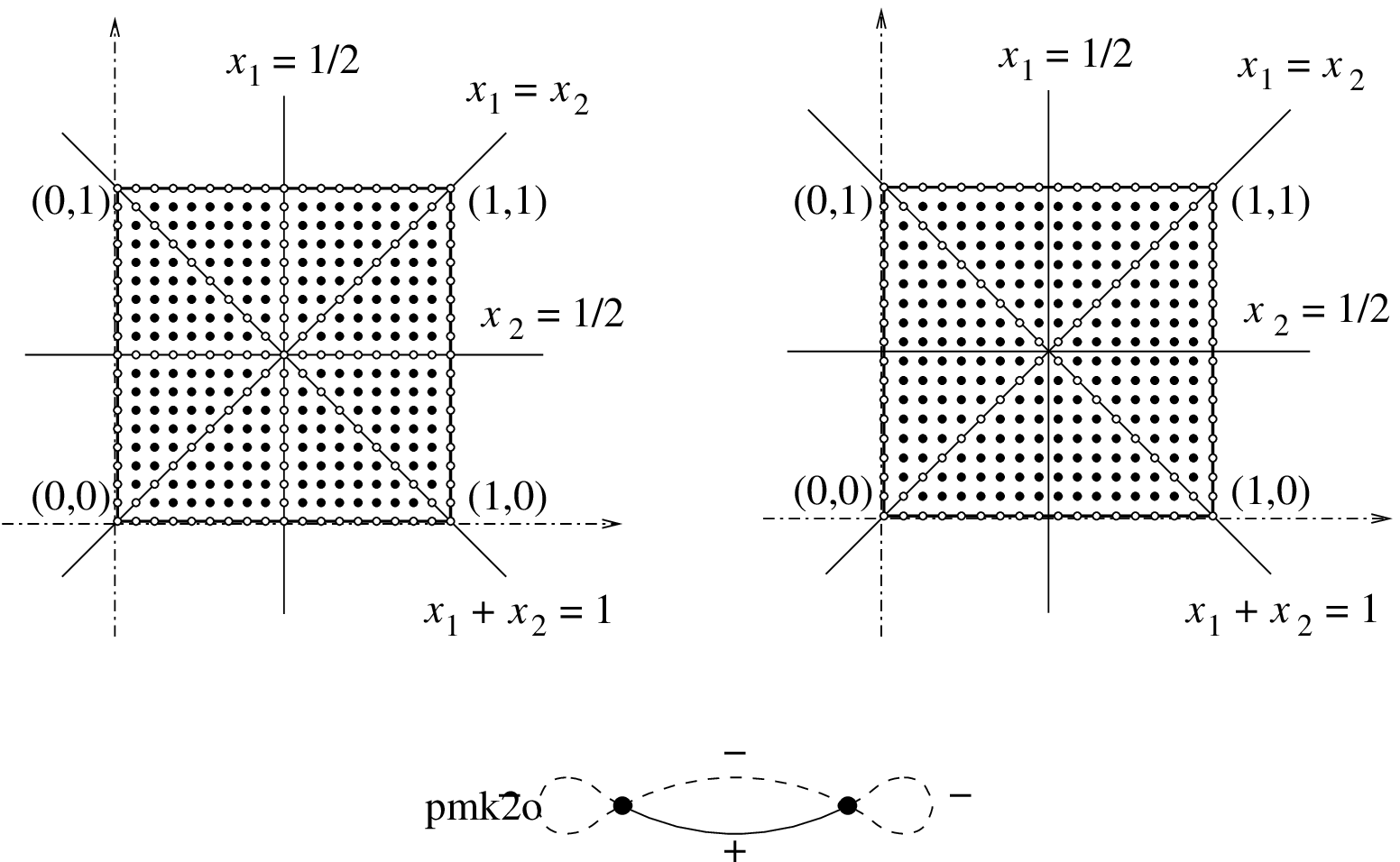}
\end{center}
\caption{Illustrating for $k=8$:  On the left, the $\frac{1}{2(k+1)}$-lattice points in $[0,1]^2$ that $k$-color the signed graph $\pm K_2^\circ$, with shifted hyperplanes.
On the right, the $\frac{1}{(2k+1)}$-lattice points that $k$-color it without $0$.} \label{F:sgcol2shift}
\end{figure}

In order to see how this is so, we must realign and rescale the whole picture so that the fundamental polytope is $P = [0,1]^n$ (just as with unsigned graph coloring) and the hyperplanes center on the point $\frac{1}{2}\bj$, where $\bj := (1,1,\hdots,1)$.  
(See Figure \ref{F:sgcol2shift}.)  
 We replace $\cH[\Sigma]$ by its translate $\cH''[\Sigma] =\cH[\Sigma] + \frac{1}{2} \bj$; that is, we add $(\frac{1}{2},\frac{1}{2},\hdots,\frac{1}{2})$ to every hyperplane.

\begin{thm} \label{T:sg}  
Let $\Sigma$ be a signed graph and let $P = [0,1]^n$.  The Ehrhart quasipolynomial of $(P,\cH''[\Sigma])$ satisfies
\begin{equation*}
(-1)^n E_{P,\cH''[\Sigma]}(-t) = E^\circ_{P^\circ,\cH''[\Sigma]}(t) = 
\begin{cases} 
\chi_\Sigma(t-1) &\text{if $t$ is even,} \\
\chi^*_\Sigma(t-1) &\text{if $t$ is odd.}
\end{cases}
\end{equation*}
\end{thm}

\begin{proof}  
An easy way to see the correctness of the expression for $E^\circ$ in terms of $\chi_\Sigma$ and $\chi^*_\Sigma$ is to translate the center of $P$ to the origin and dilate by $t$.  Then $P$ becomes $\tilde P = [-\frac{t}{2},\frac{t}{2}]^n$ and $\cH''[\Sigma]$ becomes $\cH[\Sigma]$.  What happens to $t\inv\bbZ^n$ depends on the parity of $t$.  If $t$ is even, $t\inv\bbZ^n$ becomes $\bbZ^n$ and, much as in the introduction and the proof of Theorem \ref{T:g}, $E^\circ_{P^\circ,\cH[\Sigma]}(t) =\chi_\Sigma(2c+1)$ with $c = \frac{t}{2}-1$.  When $t$ is odd, $t\inv\bbZ^n$ is transformed to $\bbZ^n+\frac{1}{2}\bj$, in which no vector has an integral entry; the number of points of this lattice in $\tilde P^\circ$ and not in $\bigcup\cH[\Sigma]$ equals $\chi^*_\Sigma(t)$ if we regard the latter as counting colorings with color set $\frac{1}{2}\{\pm 1,\pm 3,\hdots,\pm(t-2)\}$, which is an acceptable color set because it consists of $\frac12(t-1)$ colors, each with both signs, and does not contain 0.  
\end{proof} 

The effect on the geometry of the parity of $t$ is what prevents us from conveniently stating the entire result in terms of $\cH[\Sigma]$.  See Theorem \ref{T:sgcharacteristic}.

\begin{cor}[Zaslavsky {\cite[Theorem 2.2]{SGC}}]  \label{C:sgpoly} 
For a signed graph with no positive loops or loose edges, $\chi_\Sigma$ and $\chi^*_\Sigma$ are monic polynomials of degree $n$. If $\Sigma$ has a positive loop or a loose edge, $\chi_\Sigma = \chi^*_\Sigma = 0$.
\end{cor}

\begin{proof}  
The leading terms are $1x^n$ because $P$ is $n$-dimensional with volume 1.  Polynomiality is a consequence of the next lemma, by which $(P,\cH''[\Sigma])$ has denominator 1 or 2.
\end{proof}

\begin{lem}  \label{L:sgdenominator}  
If $\Sigma$ is a signed graph, $(P,\cH''[\Sigma])$ has half-integral vertices.
\end{lem}

\begin{proof}  
The flats of $\cH[\Sigma]$ correspond to partial signed partitions $(\pi,\sigma)$ of $V$.  (This description is based on \cite[Section 3]{CGL}, interpreted in light of \cite[Theorems 5.1(b) and 8.1]{SG}.)  A \emph{partial partition} is a partition of a subset of $V$.  A \emph{signed partition} is a partition $\pi$ along with, for each block $B$, a pair $[\sigma]= \{\sigma,-\sigma\}$ where $\sigma : B\to \{+1,-1\}$ is a signature on $B$.  In the correspondence $(\pi,\sigma)\mapsto s$, the flat $s$ has the equations $x_i = 0$ for $i \not\in \bigcup\pi$ and, for each block $B\in \pi$, $\sigma_ix_i = \sigma_jx_j$ if $i,j\in B$.  (In general not all subspaces of this form are flats of $\cH[\Sigma]$, the exception being the complete signed graph $\pm K_n\full$ \cite{SG}.)

A flat of $\cH''[\Sigma]$ therefore has the equations $x_i = \frac{1}{2}$ 
if $i \not\in \bigcup\pi$ and 
$\sigma_jx_j = \sigma_i x_i +\frac{1}{2} (\sigma_j-\sigma_i)$ if 
$i \underset\pi\sim j$.  The constant term in the latter is integral.  
A vertex of $(P,\cH''[\Sigma])$ is described by $n-k$ equations of these kinds, determining a $k$-flat $s$, and $k$ equations of the form $x_i = a_i \in \{0,1\}$; clearly, then, the vertex has half-integral coordinates.
\end{proof}

We say more about half integrality in relation to the incidence matrix in \cite{NNZ}.

There is a stronger conclusion if $\Sigma$ is \emph{balanced}, that is, it has no halfedges and no circles with negative sign product.  In that case $\Sigma$ is obtained from an all-positive graph by reversing the signs of all edges of a cutset, an operation called \emph{switching}.  (This was proved, in essence, by K\"onig \cite[Theorem X.10]{K}.  See \cite[Corollary 3.3]{SG} for more detail.)  We represent switching by a function $\eta : V \to \{+,-\}$ such that the cutset consists of all edges whose endpoints have opposite signs. When $\Sigma$ is balanced, obtained by switching $+\Gamma$ (where $\Gamma$ is the underlying graph of $\Sigma$), an edge has sign $\sigma(ij) = \eta(i)\eta(j)$ and a flat of $\cH[\Sigma]$ is specified by a partition $\pi$ of $V$ and equations $\eta(i)x_i = \eta(j)x_j$ when $i \underset\pi \sim j$.  

\begin{cor} [Zaslavsky {\cite[Section 2.1]{SGC}}]  \label{C:sgbalpoly}  
For a balanced signed graph, $\chi_\Sigma = \chi^*_\Sigma$.
\end{cor}

\begin{proof}  
$\Sigma$ is obtained through switching $+\Gamma$ by a switching function $\eta$.  The effect of $\eta$ on $\cH''[\Sigma]$ is to reverse coordinates:  $x_i \to 1-x_i$ if $\eta(i) = -$, but $x_i \to x_i$ if $\eta(i) = +$.  We apply $\eta$ to $P$ and $t\inv\bbZ^n$ in the same way so that switching does not alter the Ehrhart quasipolynomials.  Therefore, we may apply Lemma \ref{L:gdenominator} to $(P,\cH''[\Sigma])$.
\end{proof}

The switching equivalence of a balanced signed graph to an all-positive graph demonstrates that $(P,\cH''[\Sigma])$ then has integral vertices.  
Another proof is by observing that its equations are totally unimodular, that is, every subdeterminant is 0 or $\pm1$, as shown in \cite[Theorem 1]{Heller} and later in \cite[Proposition 8A.5]{SG}.  We omit the details.

Corollary \ref{C:sgbalpoly} is not the whole story. Going beyond Ehrhart theory, one can prove that $\chi_\Sigma \neq \chi^*_\Sigma$ when $\Sigma$ is unbalanced, by comparing the lattice $\Lat G(\Sigma)$ of closed subgraphs of $\Sigma$ to the semilattice $\Latb\Sigma$ of closed, balanced subgraphs \cite[Section 5]{SG}.  They are equal if and only if $\Sigma$ is balanced, and by \cite[Theorem 2.4]{SGC} $\chi_\Sigma =\chi^*_\Sigma$ if and only if they are equal.  Expressed in Ehrhartian terms: for a signed-graphic inside-out polytope the period of the Ehrhart quasipolynomial is equal to the denominator of $(P,\cH)$.

The relationship between the hyperplane arrangement and the chromatic polynomials of a signed graph is rather complicated.  
For a flat $u$ of $\cH[\Sigma]$ let $\Sigma(u)$ be the subgraph consisting of the edges whose hyperplanes contain $u$.

\begin{thm} [Zaslavsky {\cite[Theorem 2.4]{SGC}}] \label{T:sgcharacteristic}
For a signed graph $\Sigma$,
$$
\chi_\Sigma(c) = p_{\cH[\Sigma]}(c) 
$$
and 
$$
\chi^*_\Sigma(c) = \sum_{\substack{u\in \cL(\cH)\\ \Sigma(u)\text{ is balanced}}} \mu(\0,u) c^{\dim u} .
$$
\end{thm}

\begin{proof}
We apply Theorem \ref{T:ehrhyp} again.  Assume $t$ is even.  Then $\frac12\bj$, which belongs to every flat of $\cH''[\Sigma]$, is one of the coloring points, so (as one can easily see) every $E_{P^\circ\cap u}(t) = (t-1)^{\dim u}$.  Therefore, $E^\circ_{P^\circ,\cH''[\Sigma]}(t) = p_{\cH''[\Sigma]}(t-1)$.  This equals $\chi_\Sigma(t-1)$ by Theorem \ref{T:sg}.

If, however, $t$ is odd there are two kinds of flat.  Any flat that lies in a hyperplane $x_i=\frac12$ contains no coloring points at all.  One can see from the proof of Lemma \ref{L:sgdenominator} and the fact that a flat is balanced if and only if its signed partial partition is a partition (that is, $\bigcup\pi=V$) that these are precisely the flats that correspond to unbalanced subgraphs.  These flats therefore drop out of the sum in \eqref{E:oehrhyp}.  The other flats, which correspond to balanced subgraphs, behave as in the even case.
\end{proof}

The signed-graphic generalization of Stanley's theorem, Corollary \ref{C:gcompatible}, is also a consequence of Ehrhart theory.

\begin{cor}  [Zaslavsky {\cite[Theorem 3.5]{SGC}}]  \label{C:sgcompatible}  
The number of compatible pairs $(\alpha,x)$ consisting of an acyclic orientation $\alpha$ and a $c$-coloring of a signed graph $\Sigma$ is equal to $(-1)^n\chi_\Sigma(-(2c+1))$.  The number in which $x$ is zero-free equals $(-1)^n\chi^*_\Sigma(-2c)$.  
In particular, $(-1)^n \chi_\Sigma(-1) =$ the number of acyclic orientations of $\Sigma$.
\end{cor}

\begin{proof}[Sketch of Proof]  
We omit the details of proof because they are as in our proof of Stanley's theorem.  We omit the definitions because they are lengthy.  Acyclic orientations and compatible pairs are defined in \cite[Section 3]{SGC}.  Acyclic orientations are defined in \cite{OSG} and their correspondence to regions of $\cH[\Sigma]$ is proved in \cite[Theorem 4.4]{OSG}.
\end{proof}

\begin{prob} \label{Pr:chisigma1}
A combinatorial interpretation of $(-1)^n \chi^*_\Sigma(-1)$ would be a valuable contribution, since it would interpret the quasiconstant term $c_0(1)$ of the $t\equiv1$ polynomial.  
\end{prob}

\newcommand\G{\Gamma}
\newcommand\lcm{\operatorname{lcm}}

\section{In which we compose an integer into partially distinct parts} 
\label{agc}

A \emph{composition} of a positive integer $t$ is a representation of $t$ 
as an ordered sum of positive integers: $x_1+x_2+\cdots+x_n$.  Each $x_i$ 
is a \emph{part} of the composition.  The number of compositions of $t$ 
into $n$ parts, and the number of compositions into $n$ distinct parts, 
are classical combinatorial problems.  The intermediate cases, where 
the pairs that must not equal each other are specified by a 
graph $\G$ of order $n$ (we call such a composition \emph{$\G$-strict}), 
give another application of inside-out polytopes.  We define $c_\G(t)$ 
to be the number of $\G$-strict compositions of $t$ (into $n$ parts, since 
there is a variable for each vertex of $\G$).  

This is actually a kind of graph coloring, which we call \emph{affine 
coloring} because the colors are positive integers with a prescribed sum 
$t$ rather than simply belonging to the range from 1 to $t-1$.  
Otherwise, affine coloring is just like ordinary coloring.  From this 
viewpoint $c_\G(t)$ is the number of colorings in $t$ colors that are 
affine and proper.

There are also \emph{improper} affine colorings, where we allow the value 
$0$.  (Thus these are colorings in $t+1$ colors.)  The notion of a 
compatible acyclic orientation is the same as with ordinary coloring.  The 
corresponding kind of composition allows parts equal to 0; this is a 
\emph{weak composition} of $t$.

Let $\lambda(k) := \lcm(1,2,\ldots,k)$.

\begin{thm}\label{T:agc}
 The function $c_\G(t)$ is a quasipolynomial whose period divides 
$\lambda(n_1)$, where $n_1$ is the largest order of a component of $\G$.  
Furthermore, $c_\G(0)$ is the number of acyclic orientations of $\G$.  
More generally, $(-1)^{n-1} c_\G(-t)$, for $t\geq0$, is the number of 
pairs consisting of an arbitrary weak composition $(x_1,\ldots,x_n)$ of 
$t$ into $n$ parts and, for each level set $x\inv(k)$, $0\leq k \leq t$, 
an acylic orientation of the subgraph of $\G$ induced by $x\inv(k)$.
 \end{thm}

\begin{proof}
 A composition of $t$ can be considered as an integer point in the 
interior of the $t$-fold dilation of the standard simplex $s^{n-1}$ in 
$\bbR^n$ (the simplex that is the convex hull of the $n$ standard unit 
basis vectors).  It follows from Theorem \ref{T:quasi} that 
$c_\G(t)$ is a quasipolynomial.  The vertices of the corresponding 
inside-out polytope, $(s^{n-1},\cH[\G])$, have denominators that range 
from 1 to $n_1$.  The reason is that a flat $u$ of $\cH[\G]$ corresponds 
to a partition $\{B_1,\ldots,B_k\}$ of $V$ whose blocks induce connected 
subgraphs of $\G$.  The equations of such a flat are that the $x_i$ are 
constant on each block $B_j$.  Therefore, the vertices of $u \cap s^{n-1}$ 
are the points where all $x_i$ are 0 except on one block, $B_j$, on which they 
are equal and their sum is 1.  So, the nonzero $x_i=1/|B_j|$.  It follows 
that the denominator of $(s^{n-1},\cH[\G])$ is $\lambda(n_1)$.

The arguments about acyclic orientations are similar to those in Section 
\ref{graphs}.
\end{proof}

\section{In which we become antimagic}  \label{antimagic}

In an \emph{antimagic labelling} several sums are required to be unequal.  
The general antimagic picture starts with homogeneous, rational linear 
forms $f_1,\ldots,f_m \in (\bbR^d)^*$, which may for instance be the line 
sums of a \emph{covering clutter}: a nonvoid finite set $X$ of points 
together with a family $\cL$ of subsets, called \emph{lines}, of which 
none contains any other, and whose union is $X$.  We want to count integer 
points $x$, drawn from a bounded subset of $\bbZ^d$ which we take to be 
the set of integral points in $[0,t]^d$ or $(0,t)^d$, such that
 $$
 f_j(x) \neq f_k(x) \qquad \text{ if } j \neq k.
 $$
 We may or may not require that the coordinate values of a point be all distinct; thus we have \emph{strongly} or \emph{weakly} antimagic labellings of $[d]$.  
Let us therefore define, given the forms $f_1,\ldots,f_m$, the \emph{weak antimagic enumerator} 
\begin{enumerate}
\item[]
\begin{enumerate}
\item[$A^\circ(t) :=$] the number of integer points $x \in (0,t)^d$ such that all $f_j(x)$ are distinct,
\end{enumerate}
\end{enumerate}
and the \emph{strong antimagic enumerator} 
\begin{enumerate}
\item[]
\begin{enumerate}
\item[$A^*{}^\circ(t) :=$] the number of such points $x$ in which all entries $x_i$ are also distinct.
\end{enumerate}
\end{enumerate}
These are the open Ehrhart polynomials of inside-out polytopes with 
$$
P = [0,1]^d
$$ 
but with hyperplane arrangements of a new kind, as we now explain.

We want to think of the forms as a single function $f = (f_1,\ldots,f_m) : \bbR^d \to \bbR^m$.  The antimagic property is the requirement that $f(x) \notin \bigcup\cH[K_m]$ in $\bbR^m$.  
($K_m$ denotes the complete graph on $m$ nodes.)  
Let us imagine that $f$ is any linear transformation $\bbR^d \to \bbR^m$ and that in $\bbR^m$ we have a hyperplane $h$ that is the kernel of a homogeneous (or affine) linear functional $\phi$.  Then $\phi f$ is a homogeneous (or affine) linear functional on $\bbR^d$ defining a hyperplane $h^\sharp$, the \emph{pullback} of $h$.  Applying this construction to all the hyperplanes of an arrangement $\cH$ in $\bbR^m$ we get the \emph{pullback} $\cH^\sharp$ in $\bbR^d$.  Note that $\cH^\sharp$ might include the degenerate hyperplane $\bbR^d$, even if $\cH$ does not, since $h^\sharp$ is degenerate if and only if $h \supseteq \Image f$.  
The antimagic property of $x$ is now the statement that $x \in \bbR^d \setminus \bigcup\cH[K_m]^\sharp$.  The entries of $x$ are distinct if $x \notin \bigcup\cH[K_d]$.  Thus, the hyperplane arrangement for $A^\circ$ is $\cH[K_m]^\sharp$ and for $A^*{}^\circ$ it is $\cH := \cH[K_m]^\sharp \cup \cH[K_d]$.  To complete the preparation for our antimagic theorem, recall the \emph{multiplicity} of $x$ with respect to $\cH[K_m]^\sharp$ or $\cH$ from Section \ref{conv}, here written $m(x)$ or $m^*(x)$ for simplicity.  We define
\begin{enumerate}
\item[]
\begin{enumerate}
\item[$A(t) :=$] the sum of multiplicities $m(x)$ of all integer points $x \in [0,t]^d$,
\end{enumerate}
\end{enumerate}
and 
\begin{enumerate}
\item[]
\begin{enumerate}
\item[$A^*(t) :=$] the sum of multiplicities $m^*(x)$ of all integer points $x \in [0,t]^d$.
\end{enumerate}
\end{enumerate}

\begin{thm} \label{T:antimagic}
Given homogeneous rational linear forms $f_1,\ldots,f_m : \bbR^d \to \bbR$, no two equal, the antimagic enumerators $A^*(t)$, $A^*{}^\circ(t)$, $A(t)$, and $A^\circ(t)$ are monic quasipolynomials in $t$ that satisfy the reciprocity laws
$$
A^*(t) = (-1)^d A^*{}^\circ(-t) \quad \text{ and } \quad A(t) = (-1)^d A^\circ(-t).
$$
\end{thm}

\begin{proof}
From Theorem \ref{T:quasi}.  Distinctness of the forms ensures that antimagic points $x$ do exist so that the enumerators are not identically zero.
\end{proof}

\begin{prob} \label{Pr:antimagicregions}
Is there a combinatorial interpretation of the regions?  What is the intersection-lattice structure of $\cH[K_m]^\sharp$?  It seems improbable that any simple description could be given for arbitrary forms, but maybe there is one in a special case like that of antimagic graphs.
\end{prob}

The intersection lattices of $\cH[K_m]^\sharp$ and $\cH$ are implicated in the next theorem.

\begin{lem} \label{L:antitrans}
If all forms are distinct but have equal weight, then $\cH[K_m]^\sharp$ and $\cH$ are transverse to $[0,1]^d$.
\end{lem}

\begin{proof}
If all forms have equal weight, then $\frac12 \bj \in \bigcap\cH$.  Therefore, any flat of $\cH$ or $\cH[K_m]^\sharp$ intersects $P^\circ = (0,1)^d$.
\end{proof}

\begin{thm} \label{T:antimobius}
If all forms are different but have equal weight, then
\begin{align*}
A^*{}^\circ(t) &= \sum_{u \in \cL(\cH)} \mu(\0,u) E_{(0,1)^d \cap u}(t) \ , \\
A^*(t) &= \sum_{u \in \cL(\cH)} |\mu(\0,u)| E_{[0,1]^d \cap u}(t) \ ,
\end{align*}
where $\mu$ is the M\"obius function of $\cL(\cH)$, and there are similar formulas for $A^\circ$ and $A$ with $\cH[K_m]^\sharp$ replacing $\cH$.
\end{thm}

\begin{proof}
By Theorem \ref{T:ehrhyp} and Lemma \ref{L:antitrans}.
\end{proof}

The main examples are particular cases of antimagic labelling of covering 
clutters, especially ones from graphs.  See the survey \cite[Section 
5.4]{Gallian}.

\begin{exam}[Antimagic graphs] \label{X:antigraph}
 The edges of a simple graph are labelled by integers and we want the sum 
of all labels incident to a node to be different for every node.  The 
covering clutter here has for points the edges and for lines the sets of 
all edges incident to each node.  These examples are the most studied, 
normally with the standard label set $[q]$ if there are $q$ edges (no 
doubt because the existence question is otherwise trivial).  The one case 
that must be excluded because it has no antimagic labellings is the graph 
with just two nodes and one edge.  (See \cite{Boden} for a proof.  
\cite{Boden} calls our strong antimagic ``weak'' because it reserves the 
term ``strong'' for use of the standard label set $[q]$.)

One could generalize to bidirected graphs, although we are not aware of 
any such work.  In the form associated with a node, the labels on the 
edges are added if the edge is directed into the node and subtracted if 
not.  If the graph is directed these forms have weight zero so Theorem 
\ref{T:antimobius} does apply.
 \end{exam}

A dual example has also been studied.

\begin{exam}[Node antimagic] \label{X:antinode} 
Integers are assigned to the nodes and an edge receives the sum of its endpoint values; one wants every edge to have a different label.  The forms all have weight two so Theorem \ref{T:antimobius} applies.  In the literature normally the label set is the standard one, $[n]$ where there are $n$ nodes; see \cite{Sonntag} where the notion is generalized to hypergraphs.

In the bidirected graph generalization the rule for addition and subtraction is the same as in the preceding example.  For a directed graph, therefore, the differences of the endpoint labels are what should be distinct.  (If we could take the absolute values of these differences we would be close to the famous problem of graceful labelling \cite{HR}, but we do not see how to do that within our framework.)
\end{exam}

Combining the two labellings of a graph we obtain:

\begin{exam}[Total graphical antimagic] \label{X:antitotal} 
In a \emph{total labelling} both nodes and edges are labelled.  For antimagic one wants all node and edge sums to be different.  See \cite[Section 5.4]{Gallian}.

The bidirected generalization is as in the preceding examples.
\end{exam}

\begin{exam}[Squares antimagic, semi-antimagic, and antipandiagonal, hypercubes, etc.] \label{X:antisquares}
 These are just like magic, semimagic, and pandiagonal magic squares, 
except, of course, that the line sums must all be different.  In an 
\emph{antimagic square} the lines are the rows, the columns, and the two 
diagonals.  In a \emph{semi-antimagic square} we ignore the diagonals; but 
in an \emph{antipandiagonal square} we add all the wrapped diagonals. 
 There is a scattered literature on antimagic squares, triangles, 
pentagrams, etc., in which it is generally assumed that the labels are 
consecutive.  (See Swetz \cite[p.\ 130]{Swetz} on antimagic squares.  
What we call a pandiagonal antimagic square was introduced under the name 
``heterosquare'' by Duncan, according to \cite[p.\ 131]{Swetz}.)  One 
could extend these notions to affine and projective planes, $k$-nets, and 
hypercubes but we do not know of any such work.

Our results will also apply if one imposes symmetry on squares (or hypercubes).  By this we mean that the sum of a centrally symmetric pair of numbers is constant.  Our treatment of symmetric magic squares in \cite{MML} shows how one handles symmetry geometrically.
 \end{exam}

\begin{exam}[Small antimagic] \label{X:antismall}
We take a look at $2\times2$ semi-antimagic and antimagic squares.  

First, semi-antimagic; that is, we require each row and column sum to be different.  (This is the same as antimagic edge labelling of $K_{2,2}$.)  By inspecting the equations of the hyperplanes and facets, we conclude that the vertices of $(P,\cH)$ for weak antimagic are vertices of $P$; thus we expect a monic polynomial and indeed
$$
A^\circ(t) = t^4-\frac{22}{3}t^3+19t^2-\frac{62}{3}t+8 = \frac{(t-1)(t-2)(t-3)(3t-4)}{3} .
$$
The vertices for strong semi-antimagic, however, are half integral; thus we expect, and obtain, a monic quasipolynomial of period 2:
$$
A^*{}^\circ(t) = \begin{cases}
t^4-\frac{34}{3}t^3+45t^2-\frac{218}{3}t+38 \\ 
\qquad = \frac{(t-1)(t-3)(3t^2-22t+38)}{3} &\text{if $t$ is odd}, \\
 \\
t^4-\frac{34}{3}t^3+45t^2-\frac{218}{3}t+40 \\ 
\qquad = \frac{(t-2)(t-4)(3t^2-16t+15)}{3} &\text{if $t$ is even}.
\end{cases}
$$

Now, antimagic.  The six required inequalities imply that all entries differ.  The vertices are half integral.  The enumerators are
$$
A^\circ(t) = A^*{}^\circ(t) = \begin{cases}
t^4-12t^3+50t^2-84t+45 \\ \qquad = (t-1)(t-5)(t-3)^2 &\text{if $t$ is odd}, \\
 \\
t^4-12t^3+50t^2-84t+48 \\ \qquad = (t-2)(t-4)(t^2-6t+6) &\text{if $t$ is even}.
\end{cases}
$$
\end{exam}

To conclude we mention that the theorems apply perfectly well to limited antimagic, where only some pairs of form values need be distinct, by replacing $\cH[K_m]$ with a subgraphic arrangement $\cH[\Gamma']$ where $\Gamma' \subseteq K_m$, and to partially distinct values, $\cH[K_d]$ being replaced by $\cH[\Gamma]$ for $\Gamma \subseteq K_d$.  
Moreover, one can treat negative point values, either $|x_i|<t$ or $0<|x_i|<t$, by taking the polytope $[-1,1]^d$ and a suitable hyperplane arrangement, but the same theorems will not hold exactly since the quasipolynomials are not monic.

\section{In which subspace arrangements put in their customary appearance} \label{subspace}

An \emph{arrangement of subspaces} in $\bbR^d$ is an arbitrary finite set $\cA$ of (affine) subspaces.  (We assume all the subspaces are proper.)  We wish to generalize our results to a polytope with a subspace arrangement, along the lines taken by Blass and Sagan \cite{B-S} for graph coloring.  This is possible in part.  For instance, we can define the ``multiplicity'' of a point with respect to $\cA$, but only algebraically; it need not count anything, in fact it could be negative.

To begin with we take the situation of Section \ref{conv} in which $C$ is a bounded convex set, $D$ is a discrete set, and $\cA$ is an arrangement, now a subspace arrangement, that is transverse to $C$.  We can take over most of the definitions from Sections \ref{intro}--\ref{conv} simply by changing $\cH$ to $\cA$.  For one example, the \emph{open $D$-enumerator} of $(C,\cA)$ is 
\begin{equation*}
E^\circ_{C,\cA}(D) := \#\big(D\cap C \setm \ts\bigcup\cA\big).
\end{equation*}
There are some complications, however.  The semilattice $\cL(\cA)$, still partially ordered by reverse inclusion, is not necessarily geometric or ranked; instead it is \emph{extrinsically graded} by the rank function 
$$
\rho(u) = \codim u
$$
and the total rank $\rho(\cL) = d$, so that $u$ has extrinsic corank $\rho(\cL) - \rho(u) = \dim u$.  (The notion of extrinsic grading, without a particular name, is common in writings on subspace arrangements.)  The \emph{multiplicity} of $x\in \bbR^d$ with respect to $C$ and $\cA$ is 
\begin{equation*}
m_{C,\cA}(x) := \begin{cases} 
(-1)^d p_{\cA(x)}(-1) = 
\displaystyle{\sum_{u\in \cL(\cA):x\in u}}\mu(\0,u) (-1)^{\rho(u)} &\text{if }x\in C,\\
0 &\text{if }x\not\in C.
\end{cases}
\end{equation*}
Lemma \ref{L:mult} ensures that this agrees with the definition for hyperplane arrangements, in Section \ref{conv}.  Now we can define the \emph{closed $D$-enumerator} of $(C,\cA)$ as before:
\begin{equation*}
E_{C,\cA}(D) := \sum_{x\in D} m_{C,\cA}(x).
\end{equation*}

\begin{thm} \label{T:smobius}  
Let $C$ be a bounded, convex subset of $\bbR^d$, $\cA$ a subspace arrangement that is transverse to $C$, and $D$ a discrete set in $\bbR^d$.  Then 
\begin{equation*}
E_{C,\cA}(D) = \sum_{u\in \cL(\cA)} \mu(\0,u) (-1)^{\codim u} \#(D\cap C \cap u)
\end{equation*}
and
\begin{equation*}
E^\circ_{C,\cA}(D) = \sum_{u\in \cL(\cA)} \mu(\0,u) \#(D\cap C\cap u).
\end{equation*}
\end{thm}

\begin{proof}  
That of Theorem \ref{T:mobius}, including Lemmas \ref{L:flats} and \ref{L:faces},
goes through with obvious modifications and the understanding that an ``open face'' must be interpreted as a connected component of $u \setm \bigcup\cA^u$ but may not be simply connected, much less a cell.
\end{proof}

For the main result about subspace arrangements we adapt the notation of Section \ref{ehr}, in particular the \emph{closed} and \emph{open Ehrhart functions},
\begin{equation*}
E_{P,\cA}(t) := E_{P,\cA}(t\inv D) = \sum_{x\in t\inv D}m_{P,\cA}(x)
\end{equation*}
and
\begin{equation*}
E^\circ_{P,\cA}(t) := E^\circ_{P,\cA}(t\inv D) = \# \big( t\inv D \cap \big[P \setm \ts\bigcup\cA \big] \big).
\end{equation*}

\begin{thm}  \label{T:squasi}  
If $D$ is a discrete lattice in $\bbR^d$, $P$ is a full-dimensional $D$-fractional convex polytope, and $\cA$ is a $D$-fractional subspace arrangement, then $E_{P,\cA}(t)$ and $E^\circ_{P,\cA}(t)$ are quasipolynomials in $t$, each with period equal to a divisor of the $D$-denominator of $(P,\cA)$ and with leading term $(\vol_{P}{D}) t^d$.  We have
\begin{equation} \label{E:sreciprocity}
E^\circ_{P^\circ,\cA}(t) = (-1)^d E_{\bar P,\cA}(-t).
\end{equation}
Furthermore, if $\mathcal A$ is transverse to $P$, then 
\begin{equation}\label{E:ehrsub}
E_{P,\cA}(t) = \sum_{u\in \cL(\cA)}\mu(\0,u)(-1)^{\codim u}E_{P\cap u}(t)
\end{equation}
and
\begin{equation}\label{E:oehrsub}
E^\circ_{P,\cA}(t) = \sum_{u\in \cL(\cA)}\mu(\0,u)E_{P\cap u}(t) .
\end{equation}
\end{thm}

\begin{proof}  
The two latter equations are special cases of Theorem \ref{T:smobius}.  The reciprocity law \eqref{E:sreciprocity} follows from \eqref{E:ehrsub}, \eqref{E:oehrsub}, and standard Ehrhart reciprocity.  
\end{proof}

\begin{prob} \label{Pr:subspaceconst}
The constant term $E(0)$ does not seem to have an obvious combinatorial interpretation except in special cases, as for instance if the arrangement leaves $P$ connected, when $E(0) = \epsilon(P)$ as in ordinary Ehrhart theory.
\end{prob}

\section{In which we prove a general valuation formula}  \label{val}

A \emph{normalized valuation} on the faces of a hyperplane arrangement is a function $v$ on finite unions of open faces, with values in an abelian group, such that $v(A\cup B) + v(A\cap B) = v(A) + v(B)$ for any two such unions, or more simply 
\begin{equation*}
v(F_1\cup \cdots \cup F_k) = v(F_1) +\cdots + v(F_k)
\end{equation*}
for distinct open faces $F_1,\hdots,F_k$, and also
\begin{equation*}
v(\eset) = 0
\end{equation*}
(the normalization).  For example, if $D_0$ is a finite subset of $\bbR^d$, $v(F) = \#(D_0\cap F)$ is a valuation.  Specializing further, if $D$ is a discrete set and $C$ is a bounded convex set, then $v(F) = \#(D\cap C\cap F)$ is a valuation.  For a flat $u$ of $\cH$, set
\begin{equation*}
E_v(u) = \sum_R v(R)
\end{equation*}
summed over closed regions $R$ of $\cH^u$.  If $v(F) = \#(D\cap C\cap F)$, and $\cH$ is transverse to $C$, this is simply $E_{C \cap u,\cH^u}(D)$.  

\begin{thm} \label{T:val} 
For $u\in \cL(\cH)$ and $v$ a normalized valuation on the faces of $\cH$, 
\begin{equation*}
E_v(u) = \sum_{s\in \cL:s\geq u} |\mu(s,u)|v(s).
\end{equation*}
\end{thm}

Equation \eqref{E:mobius} is a special case.  What is different about this theorem compared to Theorem \ref{T:mobius}, besides its general statement, is the proof by M{\"o}bius inversion.  The proof is more complicated, but we think it is interesting.  

Theorem \ref{T:val} can be interpreted in terms of the M{\"o}bius algebra of $\cL(\cH)$.   The M{\"o}bius algebra $M(L)$ of a poset $L$, introduced by Solomon \cite{SMA} and developed by Greene \cite{GMA}, is the algebra (over any nice ring) generated by the elements of $L$ as orthogonal idempotents.  For $u\in L$ we define $\hat u = \sum_{s\geq u} \mu(u,s)s$.  
(Technically, this defines the  M{\"o}bius algebra of the dual poset $L^*$; but that is a difference without a difference.)  
Let $\epsilon$ denote the combinatorial Euler characteristic, $\epsilon(u) = (-1)^{\dim u}$, and let $\epsilon f$ denote the pointwise product with a function $f$.  A function defined on $\cL(\cH)$ naturally extends by linearity to the M\"obius algebra of $\cL(\cH)$.  
Theorem \ref{T:val} says that, if $v$ is a normalized valuation on $\cF(\cH)$, extended in the obvious way to $\cL(\cH)$ and then to the M\"obius algebra, then $\epsilon E_v(u) = \epsilon v(\hat u)$. 

\begin{proof}  
In effect, we use M\"obius inversion twice.  

The first time is in the \emph{semilattice of faces} of $\cH$, 
\begin{equation*}
\cF(\cH) = \{ F : F \text{ is an open face of } \cH \}
\end{equation*}
ordered by inclusion of the closures.  The maximal elements of $\cF(\cH)$ are the open regions; let $\cR(\cH)$ be the set of open regions.  We show that
\begin{equation}\label{E:valinter}
(-1)^{\dim u} v(u) = \sum_{s\in \cL:s\geq u} (-1)^{\dim s} E_v(s) \qquad \text{ for } u \in \cL(\cH).
\end{equation}
Multiplying by $(-1)^{\dim u}$, the left side equals 
\begin{equation} \label{E:lhs}
\sum_{F\in\cF(\cH^u)} v(F).
\end{equation}
The right side equals
\begin{align}
\sum_{s\geq u} &(-1)^{\dim u-\dim s} \sum_{R\in \cR(\cH^s)} v(\bar R) \notag\\
&= \sum_{s\geq u} (-1)^{\dim u-\dim s} \sum_{R\in \cR(\cH^s)} \sum_{F\in \cF(\cH^s):F\leq R} v(F) \notag\\
&= \sum_{F\in \cF(\cH^u)} (-1)^{\dim u - \dim F} v(F) \sum_{\substack{R\in \cF(\cH^u)\\ R\geq F}} (-1)^{\dim R-\dim F}.	\label{E:rsum}
\end{align}
The \emph{lattice of faces} of $\cH$, $\hat{\cF}(\cH)$, is $\cF(\cH)$ with an extra top element $\1$ adjoined.  It is known that $\hat\cF(\cH)$ is Eulerian, that is, $\mu(x,y) = (-1)^{\rk y - \rk x}$ if $x\leq y$.  Thus when $x\leq y < \1$, $\mu(x,y) = (-1)^{\dim y - \dim x}$.  The inner sum in \eqref{E:rsum} is therefore
\begin{equation*}
\sum_{R\geq F} \mu_{\hat{\cF}}(F,R) = -\mu_{\hat\cF(\cH^u)}(F,\1) = (-1)^{\dim u - \dim F},
\end{equation*}
so \eqref{E:rsum} equals \eqref{E:lhs}.

Having established \eqref{E:valinter} we invert to obtain 
\begin{equation*}
(-1)^{\dim u}E_v(u) = \sum_{s\geq u} \mu_{\cL}(u,s) (-1)^{\dim s}v(s).
\end{equation*}
Multiplying this by $(-1)^{\dim u}$ and applying Rota's sign theorem, we have the theorem.
\end{proof}

For completeness we sketch a proof that $\hat{\cF}(\cH)$ is Eulerian.  
Let $\cH_{\bbP}$ be the projectivization of $\cH$, that is, $\cH \cup \{h_\infty\}$ in $\bbP^d$ with the affine hyperplanes extended into infinity.  $\cH_\bbP$ is the projection of a homogeneous hyperplane arrangement $\cH'$ in $\bbR^{d+1}$, whose face lattice is dual to that of a zonotope, whose face lattice is Eulerian because a zonotope is a convex polytope.  
Faces of $\cH'$ other than the 0-face, $F'_0 = \bigcap \cH'$, come in opposite pairs, $F'$ and $-F'$, which project to a single face $F$ of $\cH_\bbP$.  The interval $[F_1,\1]$ in $\hat{\cF}(\cH_{\bbP})$ is isomorphic to $[F'_1,\1]$ in $\hat{\cF}(\cH')$ if $F'_1$ projects to $F_1$.  
As for $F'_0$, it projects to an infinite face.  Therefore, for any face $F$ of $\cH$, which is necessarily a finite face of $\cH_\bbP$, the interval $[F,\1]$ in $\hat{\cF}(\cH)$ is equal to $[F,\1]_{\hat{\cF}(\cH_\bbP)}$, which is isomorphic to $[F',\1]_{\hat{\cF}(\cH')}$.  
It follows that $\mu(F,\1)$ in $\hat{\cF}(\cH)$ equals $(-1)^{d+1-\dim F}$.

\section*{Acknowledgement}

We thank Jeffrey C.\ Lagarias and anonymous referees for valuable 
criticisms and suggestions about the exposition.


\end{document}